\documentclass[12pt]{amsart}
\usepackage{amsmath, amssymb, euscript, hyperref}

\newtheorem{thm}{Theorem}[section]

\newtheorem{prop}[thm]{Proposition}
\theoremstyle{definition}
\newtheorem{defn}[thm]{Definition}

\newtheorem{rem}[thm]{Remark}
\newtheorem{cor}[thm]{Corollary}

\address{Azer Akhmedov, Department of Mathematics,
North Dakota State University,
Fargo, ND, 58108, USA}
\email{azer.akhmedov@ndsu.edu}

\begin{document}

\begin{center} {\bf \Large On dense subgroups of $\mathrm{Homeo}_{+}(I)$} \end{center}

\vspace{0.6cm}

\begin{center} {\bf Azer Akhmedov} \end{center}

\vspace{0.7cm}

Abstract: {\Small We prove that a dense subgroup of $\mathrm{Homeo}_{+}(I)$ is not elementary amenable. We also show that the topological group $\mathrm{Homeo}_{+}(I)$ does not satisfy the {\em Stability of the Generators Property}, moreover, any finitely generated subgroup of $\mathrm{Homeo}_{+}(I)$ admits a faithful discrete representation in it. In the last section, we demonstrate that finitely generated dense subgroups have infinite girth.}

\vspace{1cm}

\section{Introduction}

 In this paper, we study basic questions about dense subgroups of $\mathrm{Homeo}_{+}(I)$ - the group of orientation preserving homeomorphisms of the closed interval $I = [0,1]$ -  with its natural $C_0$ metric. The paper can be viewed as a continuation of {\bf [A1]} and {\bf [A2]} both of which are devoted to the study of discrete subgroups of  $\mathrm{Diff}_{+}(I)$. 
 
 \bigskip

  Dense subgroups of connected Lie groups have been studied extensively in the past several decades; we refer the reader to {\bf [BG1], [BG2], [Co], [W1], [W2]} for some of the most recent developments. A dense subgroup of a Lie group may capture the algebraic and geometric content of the ambient group quite strongly. This capturing may not be as direct as in the case of lattices (discrete subgroups of finite covolume), but it can still lead to deep results. It is often interesting if a given Lie group contains a finitely generated dense subgroup with a certain property. For example, finding dense subgroups with property $(T)$ in the Lie group $SO(n+1,\mathbb{R}), n\geq 4$ by G.Margulis {\bf [M]} and D.Sullivan {\bf [S]}, combined with the earlier result of Rosenblatt {\bf [Ro]}, led to the brilliant solution of the Banach-Ruziewicz Problem for $\mathbb{S}^n, n\geq 4$.

 \medskip
 
 A major property that we are interested in for dense subgroups of $\mathrm{Homeo}_{+}(I)$ is not property $(T)$ but amenability (incidentally, it is not known if $\mathrm{Homeo}_{+}(I)$ has a non-trivial subgroup with property $(T)$). A very natural example of a finitely generated dense subgroup of $\mathrm{Homeo}_{+}(I)$ is R.Thompson's group $F$ in its standard representation in $\mathrm{PL}_{+}(I)$. The question about its amenability has been very popular in the last four decades. On the other hand, density of a finitely generated group in a large group  $\mathrm{Homeo}_{+}(I)$ seems to be in conflict with the amenability. We prove the following theorem.
 
 \medskip
 
 \begin{thm}\label{thm:main} An elementary amenable subgroup of $\mathrm{Homeo}_{+}(I)$ cannot be dense.
 \end{thm}
 
 \medskip
 
 Let us point out that the claim of Theorem \ref{thm:main} holds for any compact manifold. However, the case of an interval (when the manifold is connected) is in fact the hardest one; when $M\ncong I$, using density, it is straightforward to arrange a ping-pong table to show that any dense subgroup of $\mathrm{Homeo}(M)$ contains an isomorphic copy of $\mathbb{F}_2$ hence it is non-amenable. 

  \medskip 
  
  We do not know how to prove a stronger result by removing the adjective {\em elementary} from the statement of the theorem. Nevertheless, we expand on some of the ideas of the proof to obtain another fact about dense subgroups. We will state the following theorem in a much more general setting than the interval. Interestingly, in this theorem, too, the case of the interval $I$ is significantly harder; in all other cases, we again invoke a ping-pong argument (this time it is less obvious, so we present this ping-pong argument in detail). 
 
 \medskip
 
 \begin{thm}\label{thm:girth} For any compact orientable manifold $M$ with positive dimension, a finitely generated dense subgroup of $\mathrm{Homeo}_{+}(M)$ has infinite girth.
 \end{thm}
 
 \medskip
 
 As a corollary of this theorem (in the case of $M = I$), we obtain $girth(F) = \infty $ reproving the results from {\bf [AST], [Br]} and {\bf [A6]}. It follows from either of Theorem \ref{thm:main} and Theorem \ref{thm:girth} that solvable groups cannot be dense in $\mathrm{Homeo}_{+}(I)$. The theorem can be stated for an arbitrary (not necessarily orientable) compact manifold of positive dimension, by replacing $\mathrm{Homeo}_{+}(M)$ with $\mathrm{Homeo}_(M)$, with essentially the same proof.
  
 \medskip
 
 The notion of girth for a finitely generated group was first introduced in {\bf [S]} in connection with the study of Heegaard splittings of closed 3-manifolds. 
  
  \begin{defn} Let $\Gamma $ be a finitely generated group. For any finite generating set $S$ of $\Gamma $, $girth (\Gamma , S)$ will denote the minimal length of relations among the elements of $S$. Then we set $$girth (\Gamma ) = \displaystyle \mathop {\sup} _{\langle S \rangle = \Gamma , |S| < \infty }girth (\Gamma , S)$$ 
 \end{defn}
  
  \medskip
  
   Basic properties of girth have been studied in {\bf [A4]}. By the above definition, an infinite cyclic group has infinite girth, but this fact should be viewed as a degeneracy since (as remarked in {\bf [A4]}) any other group satisfying a law has a finite girth. We refer the reader to {\bf [BE], [Nak1], [Nak2], [Y]} for further studies of girth. 
 
 \medskip
 
 Let us point out that the claims of both Theorem \ref{thm:main} and Theorem \ref{thm:girth} hold for connected semi-simple real Lie groups. Indeed, it is proved in {\bf [BG1]} that any dense subgroup of a connected real semi-simple Lie group $G$ contains a non-abelian free subgroup hence it must be non-amenable. On the other hand, by the main result of {\bf [A5]}, any finitely generated linear group with a non-abelian free subgroup has infinite girth. Combining this result with Proposition 1 of {\bf [A5]} it is not difficult to show that a finitely generated subgroup of a connected real Lie group with a non-abelian free subgroup has infinite girth.  
 
 \medskip
 
 In Section 2, we discuss the so-called stability of the generator's property which also holds for simple Lie groups but we show that this property fails in $\mathrm{Homeo}_{+}(I)$.
 
 \medskip
 
 We will say that a topological group $G$ satisfies {\em Stability of the Generators Property} (SGP) if for any finitely generated dense subgroup $\Gamma $ of $G$ generated by elements $g_1, \ldots , g_n$, there exists an open non-empty neighborhood $U$ of the  identity such that if $h_i\in g_iU, 1\leq i\leq n$ then the group generated by $h_1, \ldots , h_n$ is also dense.
 
 \medskip
 
 For a topological group $G$, the SGP can be viewed as a stability of any finite generating set (in a topological sense: a subset $S\subseteq G$ generates $G$ if it generates a dense subgroup in $G$). It is immediate to see that the group $\mathbb{R}$ does not satisfy SGP. On the other hand, it is not difficult to show the SGP for connected simple Lie groups, using Margulis-Zassenhaus Lemma. This lemma (discovered by H.Zassenhaus in 1938, and later rediscovered by G.Margulis in 1968) states that in a connected Lie group $H$ there exists an open non-empty neighborhood $U$ of the identity such that any discrete subgroup generated by elements from $U$ is nilpotent (see {\bf [Ra]}). For example, if $H$ is a simple Lie group (such as $SL_2(\mathbb{R})$), and $\Gamma \leq H$ is a lattice, then $\Gamma $ cannot be generated by elements too close to the identity.  It is easy to see that (or see {\bf [A2]} otherwise) the lemma fails for $\mathrm{Homeo}_{+}(I)$. We prove the following theorem. 
 
 \medskip
   
  \begin{thm}\label{thm:sgp} The topological group $\mathrm{Homeo}_{+}(I)$ does not satisfy Stability of the Generators Property.  
 \end{thm}
 
 \medskip
 
 We indeed prove more: given any finitely generated subgroup $\Gamma $ of $\mathrm{Homeo}_{+}(I)$, and an arbitrary $\epsilon  > 0$, we show that one can find an isomorphic copy $\Gamma _1$ of $\Gamma $ generated by elements from an $\epsilon $-neighborhood of the generators of $\Gamma $ such that $\Gamma _1$ is discrete. This also shows that {\em any finitely generated subgroup of $\mathrm{Homeo}_{+}(I)$ admits a faithful discrete representation in it}. 
 
 \medskip
 
 We also prove that {\em every} finite generating set of $\mathrm{Homeo}_{+}(I)$ is indeed unstable. Furthermore, given any $n$-tuple $(g_1,\ldots , g_n)$ generating {\em a dense subgroup}, one can find another $n$-tuple $(h_1, \ldots , h_n)$ arbitrarily close to it which generates {\em a discrete subgroup}. 
 
 \bigskip      

 It is a well known fact (see {\bf [G]} or {\bf [Nav2]}) that any countable left-orderable group embeds in $\mathrm{Homeo}_{+}(I)$. We modify this argument slightly to obtain the claim of Theorem \ref{thm:sgp}.
 
 \medskip
 
  Let us emphasize that, despite the simplicity of the argument in {\bf [G]}, it does not produce a smooth embedding. Indeed, there are interesting examples of finitely generated left-orderable groups which do not embed in $\mathrm{Diff}_{+}(I)$ {\bf [Be], [Nav1]}. For the group $\mathrm{Diff}_{+}(I)$, we do not know if the property SGP holds in either $C_1$ or $C_0$ metric; it is also unknown to us if every finitely generated subgroup $\Gamma \leq \mathrm{Diff}_{+}(I)$ admits a faithful $C_1$-discrete representation in $\mathrm{Diff}_{+}(I)$. Much worse, we even do not know if $\mathrm{Diff}_{+}(I)$ contains any finitely generated $C_1$-dense subgroup at all! 
 
 \bigskip
 
 {\em Acknowledgment:} The question about the Stability of the Generators Property was brought to my attention by Viorel Nitica. It is a pleasure to thank him for a motivating conversation. 
  
\section{Instability of the Generators}

\bigskip

 In this section, we will prove Theorem \ref{thm:sgp}. For $f\in C[0,1]$, $||f||$ will denote the usual $C_0$-norm, i.e. $||f|| = \displaystyle \max _{x\in [0,1]}|f(x)|$.
 
 \medskip
 
 First, we need the notion of a $C_0$-strongly discrete subgroup from {\bf [A1]}:
 
 \medskip
 
  \begin{defn} A subgroup $\Gamma $ is $C_0$-strongly discrete if there exists $\delta > 0$ and $x_0\in (0,1)$ such that $|g(x_0)-x_0| > \delta $ for all $g\in \Gamma \backslash \{1\}$.
\end{defn} 
 
 Notice that $C_0$-strongly discrete subgroups are $C_0$-discrete. The following theorem  is stronger than Theorem \ref{thm:sgp}. 
 
 \medskip

 \begin{thm} Let $\Gamma $ be a subgroup of $\mathrm{Homeo}_{+}(I)$ generated by finitely many homeomorphisms $f_1, \ldots , f_s$, and $\epsilon > 0$. Then there exist $g_1, \ldots , g_s\in \mathrm{Homeo}_{+}(I)$ such that $\displaystyle \max_{1\leq i\leq s}||g_i-f_i|| < \epsilon $, moreover, the subgroup $\Gamma _1$ generated by $g_1, \ldots , g_s$ is $C_0$-strongly discrete, and $\Gamma _1 $ is isomorphic to $\Gamma $. 
 \end{thm} 

 \medskip
 
 {\bf Proof.} Let $(x_0, x_1, \ldots )$ be a countable dense sequence in $(0,1)$ where $x_0 = \frac{1}{2}$, and let $\delta = \frac{1}{10}\min \{\epsilon , 1\}$. Since $\Gamma $ is finitely generated it is countable and left-orderable with a left order $\prec $ such that for all $h_1, h_2\in \Gamma $, we have $h_1\prec h_2$ iff for some $n\geq 0, h_1(x_n) < h_2(x_n)$ and $h_1(x_i) = h_2(x_i)$ for all $i < n$. 
 
 \medskip
 
 Let $\gamma _0, \gamma _1, \gamma _2 , \ldots $ be all elements of $\Gamma $ where $\gamma _0 = 1$. We will build homeomorphisms $\eta _0, \eta _1, \eta _2 , \ldots $ such that they generate a subgroup $\Gamma _1$ satisfying the following conditions:
 
 \medskip
 
 (i) $\eta _0 = 1$;
 
 \medskip
 
 (ii) there exists an isomorphism $\phi : \Gamma \rightarrow \Gamma _1$ such that $\phi (\gamma _n) = \eta _n$ for all $n\geq 0$.
 
 \medskip
 
 (iii) $d_0(\gamma _n, \eta _n) < \epsilon $ for all $n\geq 1$.    

 \medskip
  
 First, we define $\eta _n(x_0)$ inductively for all $n\geq 1$. We let $\eta _1(x_0)$ to be any number in $(0,1)$ such that 
 
 \medskip
 
 (i) $\eta _1(x_0)\notin (x_0 - \delta , x_0 + \delta )$;  
 
 \medskip
 
 (ii) $|\eta _1(x_0)-\gamma _1(x_0)| < \frac{\epsilon }{2}$; 
 
 \medskip
 
 (iii) $(\eta _1(x_0) - x_0)(\gamma _1(x_0) - x_0)\geq 0$, i.e. $\eta _1(x_0)$ and  $\gamma _1(x_0)$ are on the same side of $x_0$.
 
 \medskip
 
 Now suppose $\eta _1(x_0), \ldots , \eta _n(x_0)$ are defined. To define $\eta _{n+1}(x_0)$ we consider the following three cases.
 
 \medskip
 
 {\em Case 1:} $\gamma _i \prec \gamma _{n+1}$ for all $0\leq i\leq n$.
 
 \medskip
 
  Then we let $\eta _{n+1}(x_0)$ be any number in $(\frac{1}{2}+\delta ,1)$ such that $\eta _{n+1}(x_0) > \eta _i(x_0), 1\leq i\leq n$ and $|\eta _{n+1}(x_0)-\gamma _{n+1}(x_0)| < \frac{\epsilon }{2}$ 
  
  \medskip
  
  {\em Case 2:} $\gamma _i \prec \gamma _{n+1}\prec \gamma _j$ for some $i,j \in \{1,\ldots , n\}$ where for all $k\in \{1,\ldots , n\}\backslash \{i,j\}$ either $\gamma _k \prec \gamma _i$ and $\gamma _j \prec \gamma _k$.
 
 \medskip
  
  In this case, we let $\eta _{n+1}(x_0)\in (\eta _i(x_0), \eta _j(x_0))$ and $|\eta _{n+1}(x_0) - \gamma _{n+1}(x_0)| < \frac{\epsilon }{2}$.
  
  \medskip
  
    {\em Case 3:} $\gamma _{n+1} \prec \gamma _i$ for all $0\leq i\leq n$.
 
 \medskip
 
  Then we let $\eta _{n+1}(x_0)$ be any number in $(0, \frac{1}{2}-\delta )$ such that $\eta _{n+1}(x_0) < \eta _i(x_0), 1\leq i\leq n$ and $|\eta _{n+1}(x_0)-\gamma _{n+1}(x_0)| < \frac{\epsilon }{2}$.
  
  \medskip
  
  Thus we have defined the orbit $O(x_0) = \{\eta _n(x_0) \ | \ n\geq 0\}$, and $\eta _n(x_0) \notin (x_0 - \delta , x_0 + \delta)$ for all $n\geq 1$. Then we can extend the definition of $\eta _n, n\geq 1$ to the whole $O(x_0)$ by setting $\eta _n(\eta _m(x_0)) = (\eta _n\eta _m(x_0))$ for all $m, n\geq 0$. 
  
  \medskip
  
  Now we extend the definitions of $\eta _n, n\geq 1$ to the set of all accumulation points of $O(x_0)$: let $z$ be an accumulation point of $O(x_0)$, so $z = \displaystyle \lim _{k\to \infty }z_k$ where $z_k =\eta _{n_k}(x_0), k\geq 1$. For all $n\geq 1$, we let $\eta _n(z) = \displaystyle \lim _{k\to \infty } \eta (z_k)$. 
  
  \medskip
  
  Since the set $[0,1]\backslash \overline{O(x_0)}$ is open, it is a union of countably many disjoint open intervals. Then we can extend the definition of the maps $\eta _n, n\geq 1$ affinely to the whole $[0,1]$.  
  
  \medskip
  
  By construction, the group $\Gamma _1 = \{\eta _0, \eta _1, \ldots \}$ is isomorphic to $\Gamma $, moreover, $\eta _n(x_0) \notin (x_0 - \delta , x_0 + \delta)$ for all $n\geq 1$. Thus $\Gamma _1$ is $C_0$-strongly discrete. $\square $ 
  
  \medskip
  
  \begin{cor} Any finitely generated subgroup of $\mathrm{Homeo}_{+}(I)$ admits a discrete embedding in it.
  \end{cor}  
  
  We do not know if the claim of the corollary holds for $\mathrm{Diff}_{+}(I)$ in $C_1$ metric. It is worth mentioning that not every finitely subgroup of a Lie group admits a discrete embedding in it: the group $\mathbb{Z}\wr \mathbb{Z}$ embeds in $GL(2,\mathbb{R})$ but does not embed discretely in any connected real Lie group. 
  
  \vspace{1cm}
  
  \section{Elementary amenable subgroups of  $\mathrm{Homeo}_{+}(I)$}
  
  \bigskip
  
   In this section, we prove Theorem \ref{thm:main}. First, we give a separate proof for solvable groups. The following proposition seems interesting independently.
  
  \medskip
  
   \begin{prop}\label{prop:normal} Let $\Gamma \leq \mathrm{Homeo}_{+}(I)$ be a dense subgroup, and $N$ be a non-trivial normal subgroup of $\Gamma $. Then $N$ is dense.
  \end{prop}
  
  {\bf Proof.} Let $\epsilon > 0$ and $\phi \in \mathrm{Homeo}_{+}(I)$. We can choose a natural number $n$ and $a_1, \ldots , a_n, b_1, \ldots , b_n\in (0,1)$ such that $\frac{1}{n} < \frac{\epsilon }{2}, a_i = \frac{i}{n+1}, 0\leq i\leq n+1, 0 = b_0 < b_1 < b_2 < \ldots < b_n < b_{n+1} = 1$, and the following three conditions hold:
  
  \medskip
  
  (c1) $a_i\neq b_j$ for all $i, j\in \{1, \ldots , n\}$
  
  \medskip
  
  (c2) $|b_i-\phi (a_i)| < \frac{\epsilon }{8}, 1\leq i\leq n$.
  
   \medskip
   
   (c3) $|b_{i+1}-b_i| < \frac{\epsilon }{8}, 0\leq i\leq n$.
   
  \medskip
  
  Let also $p = \mathrm{min}\{a_1, b_1\}, q = \mathrm{max}\{a_n, b_n\}$. Since $\Gamma $ is dense, it does not have a global fixed point in $(0,1)$. Then, there exists $f\in N$ such that $f(\frac{p}{2}) >  q + \frac{1-q}{2}$. 
  
  \medskip
  
  Let $c_k = f(a_k), d_k = f(b_k), 1\leq i\leq n$ and $c_0 = d_0 = q, c_{n+1} = d_{n+1} = 1$. Then $q = c_0 < c_1 < \ldots < c_n < 1$, and $q = d_0 < d_1 < \ldots < d_n < 1$. 
  
  \medskip

  Let $\delta _0 = \frac{1}{16}\mathrm{min}\{\delta _1, \delta _2, \delta _3, \epsilon \}$ where $$\delta _1 = \displaystyle \min _{0\leq i\leq n}|c_{i+1} - c_i|, \ \delta _2 = \displaystyle \min _{0\leq i, j\leq n}|a_i-b_j|, \ \delta _3 = \min \{p, 1-q\}$$ 
  
  Then there exists a positive  $\delta  < \delta _0$ such that for all $k\in \{1, \ldots , n\}$, we have $f^{-1}(I_k) \subset (b_k - \delta _0, b_k + \delta _0)$ where $I_k = (d_k - \delta ,d_k + \delta ), 1\leq k\leq n$. 
  
  \medskip

  Now, let $J_k = (b_k - \delta _0, b_k + \delta _0), L_k = (b_k - 2\delta _0 , b_k + 2\delta _0 ), 1\leq k\leq n$.
  
  \medskip
  
  Notice that $J_k$ is a subinterval of $L_k, 1\leq k\leq n$, and the intervals $L_1, \ldots , L_n,  I_1, \ldots , I_n$ are mutually disjoint. Moreover, all of the intervals $I_1, \ldots , I_n$ lie on the right side of $q + \frac{1-q}{2}$ while all of the intervals $L_1, \ldots , L_n$ lie on the left side of $q + \frac{1-q}{2}$. 
  
  \medskip 
  
  By the density of $\Gamma $, we can find $g\in \Gamma $ such that for all $k\in \{1, \ldots , n\}$ the following conditions hold:
  
  \medskip
  
  (i) $g(c_k) \in I_k$;
  
  \medskip
  
  (ii) $g^{-1}(J_k) \subset L_k$.

  \medskip
  
  Then $g^{-1}f^{-1}gf(a_k) \in L_k$ for all $k\in \{1, \ldots , n\}$. Then using conditions (c2) and (c3) we easily obtain that $||\phi - g^{-1}f^{-1}gf|| < \epsilon $. $\square $
  
  \bigskip
  
  We now observe an important corollary of Proposition \ref{prop:normal}.
  
  \begin{cor}\label{cor:solvable} A solvable subgroup of $\mathrm{Homeo}_{+}(I)$ is not dense.
  \end{cor}
  
  {\bf Proof.} Indeed, if $\Gamma $ is a solvable dense subgroup of $\mathrm{Homeo}_{+}(I)$ then it has a non-trivial normal Abelian subgroup $N$. By Proposition \ref{prop:normal}, $N$ is dense. But since $N$ is Abelian, it has a non-trivial cyclic normal subgroup $C$. Again, by Proposition \ref{prop:normal}, $C$ is dense. However, a cyclic subgroup cannot be dense. Contradiction. $\square $
  
  \bigskip
  
  \begin{rem} In fact, by an argument similar to the one in the proof of Proposition \ref{prop:normal}, one can show that an Abelian subgroup $G$ of $\mathrm{Homeo}_{+}(I)$ cannot be $\frac{1}{4}$-dense, i.e. one can find $f\in \mathrm{Homeo}_{+}(I)$ s.t. $f$ lies in a distance $\frac{1}{4}$ apart from $G$, and a solvable subgroup cannot be $\frac{1}{8}$-dense. 
  \end{rem}
  
  \medskip
  
  Our goal is now to extend the corollary to show that a dense subgroup of $\mathrm{Homeo}_{+}(I)$ cannot be elementary amenable. 
  
  \medskip
  
  For the convenience of the reader, let us recall that the class of amenable groups is closed under the following four natural processes of forming new groups out of the old ones: {\bf (I)} subgroups, {\bf (II)} quotients, {\bf (III)} extensions, and {\bf (IV)} direct unions. Following C.Chou {\bf [C]}, let us denote the class of Abelian groups and finite groups by EG$_0$. Assume that $\alpha > 0$ is an ordinal and we have defined EG$_{\beta }$ for all ordinals $\beta < \alpha $. Then if $\alpha $ is a limit ordinal, set EG$_{\alpha } = \displaystyle \mathop {\sqcup} _{\beta < \alpha }$EG$_{\beta }$ and if $\alpha $ is not a limit ordinal, set  EG$_{\alpha }$ is the class of groups which can be obtained from groups in EG$_{\alpha -1}$ by either applying process {\bf (III)} or process {\bf (IV)} once and only once. It is proved that each class EG$_{\alpha }$ is closed under processes {\bf (I)} and {\bf (II)} and EG = $\cup \{$EG$_{\alpha } : \alpha $ is an ordinal $\}$ is the smallest class of groups which contain all finite and Abelian groups and is closed under the processes {\bf (III)} and {\bf (IV)}. A group from the class EG is called {\em elementary amenable group}. Some basic and interesting properties of these groups have been studied in {\bf [C]}.
  
  \medskip
  
  A subgroup of $\mathrm{Homeo}_{+}(I)$ from class EG$_0$ is Abelian and Abelian groups are not dense by Corollary \ref{cor:solvable}. Using this fact as a base of a transfinite induction, one would want to establish the step  of it to prove that an elementary amenable subgroup is not dense. Assume that we can prove this claim for the groups of classes EG$_{\beta }$ for all $\beta < \alpha $. If $\alpha $ is a limit ordinal then by definition of EG$_{\alpha }$, any group $\Gamma $ from it belongs to a class  EG$_{\beta }$ for some $\beta < \alpha $ thus we conclude by the inductive assumption that $\Gamma $ is thin. If $\alpha $ is not a limit ordinal then there are two ways to obtain $\Gamma $ from EG$_{\alpha -1}$: {\bf (i)} $\Gamma $ is an extension of $A$ by $B$ where $A, B$ are non-trivial subgroups from EG$_{\alpha -1}$, and {\bf (ii)} $\Gamma $ is a direct union of $\{\Gamma _{\tau }\}, \Gamma _{\tau }\in $EG$_{\alpha -1}$.
   
   \medskip
   
   In Case (i), if $\Gamma $ is dense then, by Proposition \ref{prop:normal}, the non-trivial normal subgroup $B$ is also dense; but this contradicts the inductive assumption. However, in Case (ii), we are unable to carry out the step, for the following reason: a directed union of countably many nowhere dense subgroups of $\mathrm{Homeo}_{+}(I)$ can indeed be dense!  
   
   \medskip
  
  To overcome this difficulty, we would like to introduce a concept of {\em thin groups} which helps us to take care of the problem. For an integer $n$, let  \begin{displaymath}  \mathrm{sgn}(n) =  \left\{\begin{array}{lcr} 1 &  \mathrm{if} \ n > 0 \\ 0 &  \mathrm{if} \ n =0 \\  -1 &  \mathrm{if} \ n < 0    \end{array}   \right. \end{displaymath}    
  
  \medskip
  
  \begin{defn} Let $N\geq 1$ be an integer. A group $\Gamma $ is called $N$-thin if for all $a, b\in \Gamma $ there exists a word $W(a,b) = a^{n_1}b^{n_2}\ldots a^{n_{2k-1}}b^{n_{2k}}a^{n_{2k+1}}$ such that $W(a,b) = 1\in \Gamma $ where $n_2, \ldots , n_{2k}$ are non-zero, moreover, $\mathrm{sgn}(n_1) + \ldots + \mathrm{sgn}(n_{2k+1}) = 0$, and  $|\mathrm{sgn}(n_1) + \ldots + \mathrm{sgn}(n_{i})|\leq N$ for all $i\in \{1, \ldots , 2k+1\}$.
  \end{defn}
  
  \medskip
  
  In the above definition, the quantity $\displaystyle \mathop{\max}_{1\leq i\leq 2k+1}|\mathrm{sgn}(n_1) + \ldots + \mathrm{sgn}(n_{i})|$ will be called {\em the width} of the word $W(a,b)$, and the quantity $\displaystyle \mathop{\max}_{1\leq i\leq 2k+1}|n_{i}|$ will be called {\em the height} of the word $W(a,b)$
  
  \medskip
  
  \begin{defn} A group is called thin if it is $N$-thin for some $N\geq 1$.
  \end{defn}
  
  \medskip
  
  Let us observe the following important facts.
  
  \medskip
  
  \begin{prop}\label{prop:thin} (i) a subgroup of an $N$-thin group is $N$-thin;
  
  (ii) a quotient of an $N$-thin group is $N$-thin;
  
  (iii) an extension of an $N$-thin group by an $M$-thin group is $(M+N)$-thin;
  
  (iv) a directed union of $N$-thin groups is $N$-thin;
  
  (v) Abelian groups are 1-thin;
  
  (vi) finite groups are 1-thin. \ $\square $
  \end{prop}
  
  \medskip
  
   Thin groups turn out interesting from a pure group-theoretical point of view. Despite Proposition \ref{prop:thin}, not all elementary amenable groups are thin. Conversely, the class of thin groups includes interesting groups which are not elementary amenable (and not even amenable). Still, thin groups are useful in understanding the concept of amenability; here, we will limit ourselves to pointing out the following basic property of these groups.  
  
  \medskip
   
  \begin{prop}\label{prop:thindense} A thin subgroup $\Gamma $ of $\mathrm{Homeo}_{+}(I)$ is not dense.
  \end{prop} 
  
  \medskip
  
  {\bf Proof.} Let $n = 2N+2, 0 < a_1 < b_1 < a_2 < b_2 < \ldots < a_n < b_n < 1$, and $x_0\in (a_{N+1}, b_{N+1})$. Let also $$S = \{0, x_0, 1\}\sqcup \{a_1, \ldots , a_n\}\sqcup \{b_1, \ldots , b_n\}$$  $$ \mathrm{and} \ \epsilon = \frac{1}{10}\min \{x-y \ | \ x,y\in S, x\neq y\}$$   
  
  \medskip
  
  We choose two homeomorphisms $f, g\in \mathrm{Homeo}_{+}(I)$ such that the following conditions hold:
  
  \medskip
  
  (i) $Fix(f) = \{a_1, \ldots , a_n\}\cup \{0,1\}, Fix(g) = \{b_1, \ldots , b_n\}\cup \{0,1\}$;
  
  \medskip
  
  (ii) for all $x\in I, f(x) \geq x$ and $g(x)\geq x$;
  
  \medskip
  
  (iii) for all $x\in I$ if $\min \{x-z \ | \ z\in \{0,1\}\sqcup \{a_1, \ldots , a_n\}\}  > 2\epsilon $ then $f(x)-x > \epsilon $;
  
  \medskip
  
  (iv) for all $x\in I$ if $\min \{x-z \ | \ z\in \{0,1\}\sqcup \{b_1, \ldots , b_n\}\}  > 2\epsilon $ then $g(x)-x > \epsilon $.
  
  \medskip
  
  Since $\Gamma $ is dense, we can choose $\phi , \psi \in \mathrm{Homeo}_{+}(I)$ such that $||\phi - f||< \epsilon $ and $||\psi - g||< \epsilon $. Then the following conditions hold:
  
  \medskip
  
  (i) $Fix(\phi )\backslash \{0,1\} \subset \displaystyle \mathop {\sqcup} _{1\leq i\leq n}(a_i-2\epsilon ,a_i+2\epsilon )$;
  
  \medskip
  
  (ii) $\phi (x) > x$ for all $x\notin \displaystyle \mathop {\sqcup} _{1\leq i\leq n}(a_i-2\epsilon ,a_i+2\epsilon )\sqcup [0, 2\epsilon )\sqcup (1-2\epsilon , 1]$ 
  
  \medskip
  
  (ii) $Fix(\psi )\backslash \{0,1\} \subset \displaystyle \mathop {\sqcup} _{1\leq i\leq n}(b_i-2\epsilon ,b_i+2\epsilon )$;
  
  \medskip
  
  (iv) $\psi (x) > x$ for all $x\notin \displaystyle \mathop {\sqcup} _{1\leq i\leq n}(b_i-2\epsilon ,b_i+2\epsilon )\sqcup [0, 2\epsilon )\sqcup (1-2\epsilon , 1]$ 
  
  \medskip
  
  Notice that the intervals $[0, 2\epsilon ], [a_1-2\epsilon ,a_1+2\epsilon ], [b_1-2\epsilon ,b_1+2\epsilon ],  \ldots , [a_n-2\epsilon ,a_n+2\epsilon ], [b_n-2\epsilon ,b_n+2\epsilon ], [1 - 2\epsilon , 1]$ are mutually disjoint, and $x_0$ does not belong to any of them. Moreover, for a sufficiently big positive integer $m$, and for all $i\in \{2, \ldots , n-1\}$, we have $$\phi ^{-m}([b_i-2\epsilon ,b_i+2\epsilon ]\subset (a_i ,a_i+2\epsilon ), \phi ^m([b_i-2\epsilon ,b_i+2\epsilon ]\subset (a_{i+1}-2\epsilon ,a_{i+1})$$ and  $$\psi ^{-m}([a_i-2\epsilon ,a_i+2\epsilon ]\subset (b_{i-1} , b_{i-1}+2\epsilon ), \psi ^{m}([a_i-2\epsilon ,a_i+2\epsilon ]\subset (b_{i}-2\epsilon ,b_{i})$$ 
  
  \medskip
  
  We also have $$\phi ^{-m}(x_0) \in (a_{N+1}, a_{N+1}+2\epsilon ),  \phi ^{m}(x_0) \in (a_{N+2}-2\epsilon , a_{N+2}),$$ and  $$\psi ^{-m}(x_0) \in (b_{N}, b_{N}+2\epsilon ), \psi ^{m}(x_0) \in (b_{N+1} - 2\epsilon , b_{N+1})$$ 
  
  \medskip
  
  Then we let $a = \phi ^m, b = \psi ^m$, and observe that for sufficiently big $m$, $W(f,g)(x_0) \in \displaystyle \mathop {\sqcup} _{1\leq i\leq n}(a_i-2\epsilon ,a_i+2\epsilon )\sqcup \displaystyle \mathop {\sqcup}_{1\leq 1\leq n}(b_i-2\epsilon , b_i+2\epsilon )$ for all reduced words $W(a, b) =  a^{n_1}b^{n_2}\ldots a^{n_{2k-1}}b^{n_{2k}}a^{n_{2k+1}}$ where $|n_1+\ldots + n_i|\leq N$ for all $i\in \{1, \ldots , 2k+1\}$. Hence $W(x_0)\neq x_0$, then $W\neq 1\in \Gamma $. $\square $ 
  
  \medskip
   
   Not every elementary amenable group is thin, thus we cannot apply Proposition \ref{prop:thindense} to prove Theorem \ref{thm:main}. We will introduce a more subtle concept related to thinness. First, we need to introduce the notions of span (for subsgroups and elements of $\mathrm{Homeo}_{+}(I)$ as well as for subsets of $I = [0,1]$) and norm (for the elements of $\mathrm{Homeo}_{+}(I)$).
   
   \medskip
   
   \begin{defn}\label{defn:span} For all $g\in \mathrm{Homeo}_{+}(I)$ we let $$Span(g) = sup \{|J| \ : \ J \ \mathrm{is \ a \  subinterval \ of} \ (0,1), Fix(g)\cap J = \emptyset \}.$$ For subgroups $G\leq \mathrm{Homeo}_{+}(I)$ we let $Span(G) = \displaystyle \mathop{sup} _{g\in G}Span(g)$. Finally, for any subset $S\subseteq (0,1)$ we let $$Span(S) = sup\{|b-a| \ : \ a,b\in S, (a,b)\cap S= \emptyset \}.$$
   \end{defn}
   
   \medskip
   
   \begin{defn}\label{defn:norm} For all $g\in \mathrm{Homeo}_{+}(I)$, we let $$N_1(g) = sup \{\frac{1}{\epsilon } \ | \ 0<\epsilon < \frac{1}{2}, g(\epsilon ) > 1-\epsilon \ \mathrm{or} \ g(1-\epsilon ) < \epsilon \},$$ \ $$N_2(g) = \frac{1}{|g(\frac{1}{2})-\frac{1}{2}|}, N(g) = \max \{N_1(g), N_2(g)\}.$$ 
   \end{defn}

  \medskip
  
  Let us clarify that in case of $g(\frac {1}{2}) = \frac {1}{2}$ we have $N(g) = \infty $.
  
  \medskip
  
  We now consider the following technical property for subgroups of $\mathrm{Homeo}_{+}(I)$. We say a subgroup $G\leq \mathrm{Homeo}_{+}(I)$ satisfies {\em property (P)} if there exists a sequence $(h_n)_{n\geq 1}$ in $G$ such that the following conditions hold:
  
  \medskip
  
  (c1) the sequence $Span(h_n)$ is increasing and $\displaystyle \mathop{\lim }_{n}Span(h_n) \geq \frac{1}{2}Span(G)$;
  
  \medskip
  
  (c2) if there exists $h\in G$ such that $Span(h) > \frac{1}{2}$ then for all $n\geq 1$, $N(h_n) \leq \max \{2N(h), 100\}$;
  
  \medskip
  
  (c3) for all $g\in G$ there exists $N \geq 1$ such that for all $k > N$ there exists a word $W(g,h_k)$ of width at most two  such that $W(x) = 1$ for all $x\in S_k$ where $S_k\subseteq (0,1)$ and $Span(S_k) < \frac{1}{2}$. 
  
  \medskip
   
  Despite the technicalities of conditions (c1)-(c3), by a transfinite induction, it is straightforward to see that all elementary amenable subgroups of $\mathrm{Homeo}_{+}(I)$ satisfy property (P). Indeed, one can check that all three of these conditions are preserved under the operations {\bf (III)} and {\bf (IV)}. Thus, it remains to prove the following 
  
  \medskip
  
  \begin{prop} A group with a property (P) is not dense in $\mathrm{Homeo}_{+}(I)$.
  \end{prop}
  
  \medskip
 
  {\bf Proof.} Let $\Gamma $ be a dense group with a property $(P)$, i.e. conditions (c1)-(c3) hold. Then there exists $h\in \Gamma $ with $Span(h) > \frac {1}{2}$. By conditions (c1) and (c2) there exists a subsequence $(h_{n_k})_{k\geq 1}$, an element $h\in \mathrm{Homeo}_{+}(I)$, a sequence of natural numbers $(m_k)_{k\geq 1}$ and the points $p_0, p, q, q_0$ such that the following conditions hold:
  
  \medskip
  
  1) $0 < p_0 \leq p < \frac{1}{2} < q \leq q_0 < 1$,   
  
  \medskip
  
  2) $|p-q| > \frac{1}{2}$,
  
  \medskip
  
  3) $h(x) > x$ for all $x\in [p_0, q_0]$,
  
  \medskip
  
  4) $\displaystyle \mathop{\lim }_{k\to \infty } h_{n_k}^{-m_k}(\frac{1}{2}) = p, \ \displaystyle \mathop{\lim }_{k\to \infty } h_{n_k}^{-2m_k}(\frac{1}{2}) = p_0$,
  
  \medskip
  
  5) $\displaystyle \mathop{\lim }_{k\to \infty } h_{n_k}^{m_k}(\frac{1}{2}) = q, \ \displaystyle \mathop{\lim }_{k\to \infty } h_{n_k}^{2m_k}(\frac{1}{2}) = q_0$. 
  
  \medskip
  
  Now, let $f\in \mathrm{Homeo}_{+}(I)$ such that $f(x) \geq x$ for all $x\in [p_0,q_0]$, and $Fix(f)\cap [p_0,q_0] = \{\frac{p+\frac{1}{2}}{2}, \frac{\frac{1}{2}+q}{2}\}$.
  
  \medskip
  
  If $g\in \Gamma $ is sufficiently close to $f$ (such an element $g$ exists by the denseness of $\Gamma $) then for sufficiently big $m$ and $k$, we have $W(g^m, h_{n_k}^{m_k})(x) \neq x$ for all $x\in [p_0, q_0]$ where $W$ is a word of width at most two, and height at most two. Contradiction. $\square $
  
  \vspace{1cm}

    \section{Dense $\Rightarrow $ Infinite Girth}
  
  In this section, we prove Theorem \ref{thm:girth}. Without loss of generality, we may assume that $M$ is connected. Interestingly, the case of $M\cong I$ seems harder than all other case, so we will treat this case separately. 
  
   \medskip 

     Let $M\ncong I$ (So $\mathrm{dim}(M)\geq 2$ or $M\cong \mathbb{S}^1$) and $\Gamma $ be a finitely generated dense subgroup of $\mathrm{Homeo}_{+}(M)$ with a finite generating set $\{\gamma _1, \ldots , \gamma _s\}$. Let $\gamma _0 = 1$. We can choose distinct points $p, a, r\in M$ and $\beta \in \mathrm{Homeo}_{+}(M)$ such that the sets $$\{\beta _i^{\epsilon }(a) : 0\leq i\leq s, \epsilon \in \{-1,1\}\}, \{\beta _i^{\epsilon }(r) : 0\leq i\leq s, \epsilon \in \{-1,1\}\}, \{p, a, r\}$$ are mutually disjoint, where $\beta _i = \beta \gamma _i, 0\leq i\leq s$; moreover, if $M \cong \mathbb{S}^1$, then the sets $$\{a\}\sqcup \{\beta _i^{\epsilon }(a) : 0\leq i\leq s, \epsilon \in \{-1,1\}\} \  \mathrm{and} \ \{r\}\sqcup \{\beta _i^{\epsilon }(r) : 0\leq i\leq s, \epsilon \in \{-1,1\}\}$$ lie in disjoint arcs and $p$ lies between these arcs.
     
     \medskip 
     
     Then we can choose $\alpha \in \mathrm{Homeo}_{+}(M)$, a natural number $N$ and disjoint open neighborhoods $U_a, U_r$ of $a, r$ respectively, such that for all $n > N$, $\alpha ^n(p)\in U_a, \alpha ^{-n}(p)\in U_r$, moreover, $\alpha ^n(\Omega )\in U_a, \alpha ^{-n}(\Omega )\in U_r$, and $(U_a\sqcup U_r)\cap \Omega = \emptyset , p\notin (U_a\sqcup U_r)\sqcup \Omega $ where $$ \Omega = \displaystyle \mathop{\cup }_{0\leq i\leq s, \epsilon \in \{-1,1\}}\beta _i^{\epsilon }(U_a\sqcup U_r) .$$

   \medskip 

    Then, for all $m\geq 2$, by taking the generating set $$S_m = \{\alpha , \alpha ^{mN}\beta _0\alpha ^{mN}, \alpha ^{2mN}\beta _1 \alpha ^{2mN}, \dots ,\alpha ^{(s+1)mN}\beta _s\alpha ^{(s+1)mN} \}$$
we observe there is no relation of length less than $m$ among the elements of $S_m$. In fact, for any word $W$ of length less than $m$ in the alphabet $S_m$, we have $W(p)\in U\sqcup V$, thus $W\neq 1\in \Gamma $.  Hence $\mathrm{girth}(\Gamma ) \geq m$. Since $m$ is arbitrary, we conclude that $\mathrm{girth} (\Gamma ) = \infty $. 
     
   \medskip

  Now, we are considering the case of $M\cong I$. Let $\Gamma $ be a finitely generated dense subgroup of $\mathrm{Homeo}_{+}(I)$, $m$ be a positive integer and $\{\gamma _1, \ldots , \gamma _s\}$ be a finite set of generators of $\Gamma $. We will find $\eta \in \Gamma $ such that the generating set $\{\eta , \eta ^m\gamma_1 \eta ^m, \ldots , \eta ^m\gamma _s\eta ^m\}$ has no relation of length less than $m$.
   
   \medskip
   
   Let $F_1$ be a free group formally generated by letters $\gamma , \gamma _1, \ldots , \gamma _s$. (by abusing the notation, we treat the elements $\gamma _1, \ldots , \gamma _s$ of $\Gamma $ also as the elements of $F_1$.) Let also $a_0 = \gamma , a_1 = \gamma ^m\gamma_1 \gamma ^m, \ldots , a_s = \gamma ^m\gamma_s \gamma ^m$, and $F_2$ be free group formally generated by $a_0, a_1, \ldots , a_s$. (so both $F_1$ and $F_2$ are free groups of rank $s+1$.) 
   
   \medskip
   
   Let $W_1, \ldots , W_N$ be all reduced words in the free group $F_2$ of length at most $m$. These words can be written as reduced words $V_1, \ldots , V_N$ in the free group $F_1$ where each word has length at most $m(2m+1)$.  
   
   \medskip
   
   Let $A = \{\gamma _1, \gamma _1^{-1} \ldots , \gamma _s, \gamma _s^{-1}\}$. We will view $A$ as a symmetrized generating set of the group $\Gamma $, and also as a finite subset of $F_1$. We will build disjoint finite subsets $S^{(0)}, S^{(1)}, \ldots , S^{(N)}$ of $(0,1)$ and define an increasing map $f:\displaystyle \mathop {\sqcup} _{i=0}^NS^{(i)}\rightarrow (0,1)$ (i.e. if $x, y \in \displaystyle \mathop {\sqcup} _{i=0}^NS^{(i)}$ and $x < y$ then $f(x) < f(y)$) inductively as follows: 
   
   \medskip
   
   First, we let $S^{(0)} = \{\frac{1}{2}\}$, and $f(\frac {1}{2})\notin \displaystyle \mathop {\sqcup}_{g\in A}g(\frac {1}{2})$. 
   
   \medskip
   
   Suppose now the subsets $S^{(0)}, S^{(1)}, \ldots , S^{(k-1)}$ are chosen and the map $f$ is defined on $\displaystyle \mathop {\sqcup} _{i=0}^{k-1}S^{(i)}$. We will describe how to define $S^{(k)}$ and extend the map $f$ to $\displaystyle \mathop {\sqcup} _{i=0}^{k}S^{(i)}$.
   
   \medskip
   
   Assume that $V_k$ has length $n$ as a reduced word in the free group $F_1$, and let $V_k = c_{n}\ldots c_2c_1$ where $c_i\in \{\gamma , \gamma ^{-1}, \gamma _1, \gamma _1^{-1}, \ldots , \gamma _s, \gamma _s^{-1}\}$, and $U_i = c_i\ldots c_2c_1, 1\leq i\leq n$, (so $U_1, \ldots , U_n$ are suffixes of $V_k$ where the reduced word $U_i$ has length $i, 1\leq i\leq n$). We define the set $S^{(k)} =\{x_0, \ldots , x_n\}$ itself and the map $f$ on it [i.e. the sequence $f(x_0), \ldots, f(x_n))]$ inductively as follows:
   
   \medskip
   
   We let $x_0$ be any point in $(0,1)$ such that $$x_0\notin \displaystyle \mathop {\sqcup} _{g\in A}g(\displaystyle \mathop {\sqcup} _{1\leq i\leq k-1}S^{(i)})\cup f(\displaystyle \mathop{\sqcup} _{1\leq i\leq k-1}S^{(i)})$$ Then we define $f(x_0) = y_0$ such that for all $x\in \displaystyle \mathop {\sqcup} _{1\leq i\leq k-1}S^{(i)}$, we have $f(x) < y_0$ iff  $x < x_0$ (so we extend the domain of $f$ such that it stays being an increasing map). 
   
   \medskip
   
   Now assume that $x_1, \ldots , x_r$ and $f(x_1), \ldots , f(x_r)$ are defined. 
   
   \medskip
   
   We consider two cases:
   
   \medskip
   
   {\em Case 1.} $U_{r+1}$ starts with $\gamma ^{\pm 1}$.
   
   \medskip
   
    Then we choose $x_{r+1}$ to be any point not in $$\displaystyle \mathop{\sqcup} _{g\in A}g(\displaystyle \mathop{\sqcup} _{1\leq i\leq k-1}S^{(i)})\cup f(\displaystyle \mathop{\sqcup} _{1\leq i\leq k-1}S^{(i)})\mathop{\sqcup} \{x_0, \ldots , x_r\}\cup \{f(x_0), \ldots , f(x_r)\}$$ and let $f(x_{r+1}) = y_{r+1}$ where for all $x\in \mathop{\sqcup} _{1\leq i\leq k-1}S_i\sqcup \{x_0, \ldots , x_r\}$, we have $f(x) < y_{r+1}$ iff $x < x_{r+1}$. 
   
   \medskip
   
   {\em Case 2.} If $U_{r+1}$ starts with some $g \in A$. 
   
   \medskip
   
   Then we let $x_{r+1} = g(x_r)$ and define $f(x_{r+1}) = y_{r+1}$ where for all $x\in \sqcup _{1\leq i\leq k-1}S_i\sqcup \{x_0, \ldots , x_r\}$, we have $f(x) < y_{r+1}$ iff $x < x_{r+1}$. 
   
   \medskip
   
   Thus we have constructed finite sets $$S^{(1)} = \{x_0^{(1)}, x_1^{(1)}, \ldots , x_{l_1}^{(1)}\}, \ldots , S^{(N)} = \{x_0^{(N)}, x_1^{(N)}, \ldots , x_{l_N}^{(N)}\}$$ corresponding to the words $V_1, \ldots , V_N$ respectively and a map $f:\displaystyle \mathop {\sqcup} _{i=1}^NS^{(i)}\rightarrow (0,1)$ such the following conditions hold:
   
   \medskip
   
   (i) $S^{(i)}$ consists of $(l_i+1)$ points where $l_i$ is the length of $V_i$ as a reduced word in $F_1$;
   
   \medskip
   
   (ii) $S^{(1)}, \ldots , S^{(N)}$ are mutually disjoint; 
   
   \medskip
   
   (iii) for all $1\leq i\leq N$, if $V_i = d_{l_i}\ldots d_2d_1$ where $d_j\in \{\gamma , \gamma ^{-1}, \gamma _1, \gamma _1^{-1}, $ \\ $\ldots , \gamma _s, \gamma _s^{-1}\}, 1\leq j\leq l_i$, then $c_j\ldots c_1(x_0^{(i)}) = x_j^{(i)}, 1\leq j\leq l_i$ where  
   
   \begin{displaymath}  c_j =  \left\{\begin{array}{lcr} d_i &  \mathrm{if} \ d_i\in \{\gamma _1, \gamma _1^{-1}, \ldots , \gamma _s, \gamma _s^{-1}\} \\ f &  \mathrm{if} \ d_i = \gamma \\  f^{-1} &  \mathrm{if} \ d_i = \gamma ^{-1}   \end{array}   \right. \end{displaymath}

   \medskip
   
   (iv) $f:\sqcup _{i=0}^NS^{(i)}\rightarrow (0,1)$ is an increasing function. 
   
   \bigskip
   
   By condition (iv), $f$ can be extended to a homeomorphism $\eta \in \mathrm{Homeo}_{+}(I)$.
   
   \bigskip

   We claim that there is no relation of length less than $m+1$ among $\eta , \gamma _1, \ldots , \gamma _s $. Indeed, let $W$ be a reduced word of length at most $m$ in the alphabet $a_1, \ldots , a_s$. Then $W$ can be written as a reduced word $V_i = V_i(\gamma , \gamma _1, \ldots , \gamma _s)$ for some $i\in \{1, \ldots , N\}$. Then $V_i(\eta , \gamma _1, \ldots , \gamma _s)$ is a homeomorphism and $V_i(\eta , \gamma _1, \ldots , \gamma _s)(x_0^{(i)}) = x_{l_i}^{(i)} \neq x_0^{(i)}$. 
   
   \medskip
   
   Let $0 < \epsilon < \mathrm{min}_{1\leq i\leq N}\frac{|x_{l_i}^{(i)}-x_1^{(i)}|}{2}$. The homeomorphism $\eta $ is not necessarily in $\Gamma $, but let us recall that $\Gamma $ is dense in $\mathrm{Homeo}_{+}(I)$. Then, if $\xi \in \Gamma $ is sufficiently close to $\eta $ we will have $V_i(\xi , \gamma _1, \ldots , \gamma _s)(x_1^{(i)}) \in (x_{l_i}^{(i)}-\epsilon , x_{l_i}^{(i)}+\epsilon )$ for all $i\in \{1, \ldots , N\}$. Hence $V_i(\xi , \gamma _1, \ldots , \gamma _s)\neq 1\in \Gamma $, for all $i\in \{1, \ldots ,N\}$. Then there is no relation of length less than $m$ among the elements of the generating set $\{\xi , \xi ^m\gamma _1\xi ^m, \ldots , \xi ^m\gamma _s\xi ^m \}$. Thus $girth(\Gamma ) \geq m$. Since $m$ is arbitrary, we conclude that $girth (\Gamma ) = \infty $. $\square $

   \vspace{1cm}
   
   {\bf R e f e r e n c e s :}
  
  \bigskip
  
  [A1] Akhmedov, A. \ On free discrete subgroups of Diff(I). \ {\em Algebraic and Geometric Topology}, \ {\bf vol.4}, (2010) 2409-2418. 
  
  \medskip
  
  [A2] Akhmedov, A. \ A weak Zassenhaus Lemma for discrete subgroups of Diff(I). {\em Algebraic and Geometric Topology,} {\bf 14}, no.1, (2014) 539-550.

  \medskip
  
  [A3] Akhmedov, A. \ Questions and remarks on discrete and dense subgroups of Diff(I). \ {\em Journal of Topology and Analysis}, {\bf vol. 6}, no. 4, (2014), 557-571 
  
  \medskip
  
  [A4] Akhmedov, A. \ On the girth of finitely generated groups. \ {\em Journal of Algebra}, {\bf 268}, (2003), no.1, 198-208.
  
  \medskip
  
  [A5] Akhmedov, A. \ The girth of groups satisfying Tits Alternative. \ {\em Journal of Algebra}, {\bf 287}, (2005), no.2. 275-282.
  
  \medskip  
  
  [A6] Akhmedov, A \ Girth alternative for subgroups of $\mathrm{PL}_{+}(I)$. \  Glasgow Mathematical Journal, {\bf 67, issue 1} (2025) 1-10. 
  
  \medskip
  
  [AST] Akhmedov, A., Stein, M., Taback, J. \ Free limits of Thompson's group $F$. Geometriae Dedicata, {\bf vol. 11} (2011) no.1, 163-176.
  
  \medskip
  
  [BE] Bartholdi, L., Erschler, A. \ Ordering the space of finitely generated groups \ http://arxiv.org/abs/1301.4669
  
  \medskip
  
   [Be] Bergman, M. G.  \ Right orderable groups that are not locally indicable. \ {\em Pacific Journal of Mathematics.} {\bf 147} Number 2 (1991), 243-248.
		
		\medskip
		
		[BG1] Breuillard, E., Gelander, T. \ On dense free subgroups of Lie groups. \ {\em Journal of Algebra}, {\bf vol.261}, no.2, (2003) 448-467.
		
		\medskip
		
		[BG2]  Breuillard, E., Gelander, T. \ A topological Tits alternative, {\em Annals of Mathematics} {\bf 166} no. 2 (2007) 427-474.

		\medskip
		
		[Br] Brin, M.: Free group of rank two is a limit of Thompson's group $F$.  Groups Geometry Dynamics, {\bf 4} (2010) no.3, 433-454.
		
		\medskip
		 
     [Ch] Chou, C. \ Elementary amenable subgroups, {\em Illinois Journal of Mathematics}, {\bf vol. 24}, issue 3 (1980), 396-407. 

  \medskip
	
	   [Co] Cornulier, Yves De. \ Dense subgroups with Property (T) in Lie groups. \ {\em Comment. Math. Helv.} {\bf 83(1)}, 55-65, 2008. 
	   
	 \medskip

	   [G] Ghys, E. \ Groups acting on the circle. \ {\em Enseign. Math. (2)} 47(3-4):329-407, 2001.
		
		\medskip
	 
	 [M] Margulis, G. \ Some remarks on invariant means, \ {\em Monatshefte f\"ur Mathematik} {\bf 90} (3): 233-235,
	 
	 \medskip
	 
	 [Nak1] Nakamura, K.  \ Ph.D Thesis, University of California Davis, 2008.
	
	\medskip
	
	[Nak2] Nakamura, K. \ The Girth Alternative for Mapping Class Groups. \ http://arxiv.org, arXiv:1105.5422.
	 	
	 	\medskip
	 	
		 [Nav1] Navas, A. \ A finitely generated locally indicable group with no faithful action by $C^1$ diffeomorphisms of the interval. {\em Geometry and Topology}, {\bf 14} (2010) 573-584.
			  
   \medskip
   
   [Nav2] Navas, A. Groups of Circle Diffeomorphisms. \ Chicago Lectures in Mathematics, 2011. {\em http://arxiv.org/pdf/0607481} 
   
   \medskip
   
   [Ra] Raghunathan, M.S. {\em Discrete subgroups of semi-simple Lie groups}. \ Springer-Verlag, New York 1972. Ergebnisse der Mathematik und ihrer Grenzgebiete, Band 68. 
  
  \medskip
    
 [Ro] Rosenblatt, J. \ Uniqueness of Invariant Means for Measure-preserving Transformations, \ {\em Transactions of AMS}, {\bf 265} (1981) 623-636.
  
  \medskip
  
  [S] Sullivan, D \ For $n > 3$ there is only one finitely additive rotationally invariant measure on the n-sphere on all Lebesgue measurable sets, \ {\em Bulletin of the American Mathematical Society}, {\bf 4} (1): 121-123
  
  \medskip
  
 [Wi1] Winkelmann, J. \ Generic subgroups of Lie groups. \ {\em Topology}, {\bf 41}, (2002), 163-181.
 
 \medskip
 
 [Wi2] Winkelmann, J. \ Dense random finitely generated subgroups of Lie groups. \  http://arxiv.org/abs/math/0309129
 
 \medskip
 
 	[Y] Yamagata, S. \ {\em The girth of convergence groups and mapping class groups.} \ Osaka Jounral of Mathematics, {\bf volume 48}, no.1, (2011).
        
\end{document}